\documentclass{article}
\usepackage[utf8]{inputenc}

\title{ZF*-Extensionality interprets full ZF}
\author{Zuhair Al-Johar }
\date{January 29 2018}

\begin{document}

\maketitle

\section{Introduction}

Dana Scott had shown that removing Extensionality from ZF set theory formalized in the customary manner would weaken it down to Zermelo set theory. The following proof is my personal attempt to solve the question of whether we can have a version of replacement that can withstand removal of Extensionality when ZF is formalized with it instead of the standard axiom schema of replacement. It is known that this can be done with axiom scheme of collection. However, arguably collection is built on an intuition that is different from that of replacement. Dana Scott had investigated that, the general lines of the proof here are similar to his; however, the form of Replacement suggested here is different.

\section {Proof}

By ZF* it is meant ZF formalized with Replacement* while keeping other axioms. The exact formulation of the axioms is: 
\smallskip

Extensionality: $\forall  X \ \forall  Y  (\forall z(z \in X \iff z \in Y) \Rightarrow  X=Y)$

Foundation:  $[\exists   X (X \in   A)]   \Rightarrow    \exists   Y \in   A (\neg \exists   Z \in   A (Z \in   Y))$

Pairing:  $\forall  A \ \forall  B \ \exists   X \ \forall  Y (Y \in   X  \iff    Y=A \lor Y=B)$

Union: $ \forall  A \ \exists   X \ \forall  Y (Y \in   X  \iff     \exists   Z \in   A (Y \in   Z))$

Power:  $\forall  A \ \exists   X \ \forall  Y (Y \in   X  \iff     \forall  Z \in   Y (Z \in   A))$

Infinity:  $\exists N ( \emptyset \in N \wedge  \forall  n \in N \exists m \in N \forall y (y \in m \iff y \in n \lor y=n))$
\smallskip

Replacement*: if $\phi(X,Z)$ is a formula in which only symbols “X”, “Z” occur free, and those only occur free, then: $$\forall  A \ \exists   B \ \forall  Y [Y \in   B  \iff     \exists   X \in   A \ \forall  Z(Z \in   Y \iff   \phi(X,Z))]$$,
is an axiom.
\medskip

CLAIM: ZF*-Ext. Interprets ZF 
\medskip

The idea of the proof is to first prove ZFA over the full domain of ZF*-Ext., then we proceed from there to prove full ZF over the domain of Pure sets of ZF*-Ext. Then we prove that regularity is dispensable with by building cumulative hierarchies, and restricting the domain to elements in them.\smallskip

The exact formulation of the version of ZFA to be interpreted is given as follows:\smallskip

Add a primitive one place predicate $``At"$ to mean ``is an atom" to the language of ZF, then define``is a set" as the $``\neg At"$, then add the following axioms:
\smallskip

Atoms:  $\forall  x (At(x)   \Rightarrow   \neg \exists   y (y  \in    x))$ [i.e. all atoms are empty]

Extensionality:  $ \forall  sets \  x,y  ( \forall  z(z  \in    x  \iff    z  \in    y)   \Rightarrow   x=y)$

Separation:  $\forall  A  \ \exists   set \ x \ \forall  y (y \in x \iff y \in A \wedge \phi(y))$; if $x$ not occurring in $\phi$

Pairing:  $\forall  A \ \forall  B \ \exists   x \ \forall  y (y \in   x  \iff    y=A \lor y=B)$

Union: $ \forall  A \ \exists    x \ \forall  y (y  \in    x  \iff     \exists   z  \in    A (y  \in    z))$

Power:  $\forall  A  \ \exists    x \ \forall  y (y  \in    x  \iff    \forall  z  \in    y (z  \in    A))$

Infinity:  $\exists N : \emptyset \in N   \wedge  \forall  n \in N \exists m \in N  \forall y (y \in m \leftrightarrow y \in n \lor y=n)$

Replacement: if $\phi(x,y)$ is a formula not mentioning $B$, then each closure of:$$\forall A [(\forall x \in A \exists! y (\phi(x,y))  \Rightarrow \exists B( set(B) \land \forall y (y \in B \iff \exists x \in A (\phi(x,y))))]$$; is an axiom. \smallskip 

Regularity:$[\exists   x (x \in   A)]   \Rightarrow    \exists   y  \in   A (\neg \exists   z \in   A (z \in   y))$ \smallskip

This allows for the existence of a Proper class of atoms.
\smallskip
 
Although this is enough to prove the claim of this article, however, some important definitions related to this proof in the style of this method would still be presented as well as the general known outline of that proof.
\smallskip

The crux of the proof is to use Marcel Crabbe approach to interpreting equality and membership that he used in interpreting NFU in SF. So equality will be interpreted as co-extensionaity, and membership as membership in $sets$, where $sets$ are defined as unions of equivalence classes under co-extensionality, denote those as 
$``=^*" ; `` \in^*"$; formally stated those are:
\smallskip

\emph{Define:} $x =^* y  \iff     \forall  z (z  \in    x  \iff    z  \in    y)$

\emph{Define:} $set(y)  \iff     \forall  m \forall n (m =^* n   \Rightarrow   [m  \in    y  \iff    n  \in    y]) $

\emph{Define:} $x  \in^* y  \iff    set(y) \wedge x  \in    y$ \smallskip

From hereafter we shall denote any object in the universe of discourse of this theory by “class” and the term “set” shall be reserved to classes meeting the above definition of set, so Ur-elements (i.e.; Atoms) are classes that are not sets.
\smallskip

Now we can see that the ``atomic" formulas $x\in^*\alpha ; \  x =^*\alpha $, are preserved under co-extensionality, i.e. if the truth value of $x \in^   *\alpha$  is $K $(for whatever $x$ and $\alpha$), then if we replace $x$ by a co-extensional object to $x$ and we replace $\alpha$ by a coextensional object to $\alpha$, then still the resulting formula would have the truth value $K$, and of course the same applies to atomic formula $x =^ * \alpha$ no doubt, this is to be called as ``truth preservation under co-extensional replacements". This entails that for any formula in the language that only uses $`` =^ *", `` \in^ *"$ as predicates, then its truth is preserved under co-extensional replacements! The reason is because the truth of any formula is a function of the truth of its atomic sub-formulas, so if the truth of latter ones is not affected by such replacements then the truth of the whole formula shall not be affected. Call this Lemma 1, formally this is:\smallskip

\textbf{Lemma 1:} if $\phi(y)$ is a formula in which $y$ is free and only occur as free, and in which $x$ does not occur, and $\phi(x)$ is the formula obtained from $\phi(y)$ by merely replacing each occurrence of $y$ by $x$, then: $$ \forall  x \forall y (x=^*y \Rightarrow [\phi(y) \iff  \phi(x)]) $$

\textbf{Corollary 1:} for any formula $\phi$ that only uses $``=^*", `` \in ^*"$ as predicates, we have: $$ \forall  x [ \forall  y (y \in    x  \iff    \phi) \Rightarrow  set(x)] $$

Proof:
Let $z=^*y$,

If $y \in    x$, then \  $\phi(y)$ \ [Definition of x]

So $\phi(z)$ \ [Lemma 1]

Then $z \in    x$ \ (Definition of x)

If $\neg y \in    x$, \ then \ $\neg \phi(y)$ \ [Definition of x]

So $\neg \phi(z)$ \ [Lemma 1]

Then \ $\neg z \in x $ \ [Definition of x]

So  $\forall  z \forall y (z=^*y   \Rightarrow   [y \in    x  \iff    z \in    x])$

QED.\smallskip

\textbf{Corollary 2:} $\forall x (set(x)  \Rightarrow    \forall m (m  \in    x  \iff    m  \in^* x))$

Proof: straightforward from definition of  $\in^*$.\smallskip

\textbf{Corollary 3:} for any formula $\phi$ that only uses $``=^*", `` \in   ^*"$ as predicates, then:$$  \forall  y (y \in    x  \iff    \phi)  \iff     \forall  y (y \in^* x  \iff    \phi)$$
 
Proof: Corollary 1, 2.\smallskip

\textbf{Proposition 1:}   $\forall  A  \exists   B  \forall  y (y \in    B  \iff     \exists   x \in    A (y=^*x))$ \smallskip

Proof: simply let $\phi(x,z)$ in Replacement* to be $z \in   x$, and B would
be a set of all $ y$s' that are coextensional with the $x$s' in A.\smallskip

\textbf{Subsidiary 1:} if $\phi$ is a formula in which $x,y$ do not occur, then: $$\forall x \forall y : x =^* y \iff [\forall z (z \in y \iff \phi) \iff  \forall z( z \in x \iff \phi)]$$; Proof: First order logic with equality and a symbol for membership $\in$ \smallskip

\textbf{Subsidiary 2:} $X \subset A \land X =^* Y \Rightarrow Y \subset A$.\smallskip

Proof: $X =^* Y \Rightarrow Y \subset X$, and since $\subset$ is transitive, then $Y \subset A$ \smallskip

\textbf{Lemma 2:} $\forall x [\forall y (y \in x \iff \forall z (z \in y \iff \phi)) \Rightarrow set(x)]$\smallskip

Proof: Let $k =^* y$

If $y \in x$, then $\forall z (z \in y \iff \phi)$ [Definition of x]

Then $\forall z (z \in k \iff \phi)$ [Subsidiary 1]

Then $k \in x$ [Definition of x]

If $y \not \in x$, then $\neg \forall z (z \in y \iff \phi)$ [Definition of x]

Then $\neg \forall z (z \in k \iff \phi)$ [Subsidiary 1]

Then $k \not \in x$ [Definition of x]

So $\forall y \forall k (k=^* y \Rightarrow [k \in x \iff y \in x])$

QED.\smallskip

 $\in$-\emph{Separation:} if $\phi$ is a formula in which $x$ does not occur, and $y$ only occur 
 
 free, then:  $\forall  A  \exists   x  \forall  y (y \in    x  \iff    y  \in    A \wedge \phi)$  is an axiom. \smallskip

Proof: simply let $\phi(x,z)$ in Replacement* be [$z=x \wedge \phi (x)$] , then any class resulting from that replacement would be a $set$ that can only contain all empty sets and all singletons of elements $x$ of $A$ that satisfy $\phi(x)$, now take any union of it and we get a class of all $x  \in    A$ satisfying $\phi$. 
\smallskip

Now clearly from $\in$-separation above we can prove the following: $$\forall A \exists x \forall y (y \in x \iff y \in^* A \land \phi)$$; where $\phi$ only uses $=^*,  \in^*$, then clearly the asserted object $x$ is a set (Corollary 1). So $\in^*$-Separation is \emph{proved}! [Corollary 3]. 
\smallskip

The proof of interpretation of Weak Extensionality, Foundation, Pairing, Union, and Power axioms of ZFA is rather straightforward (all written using predicates $=^*, \in^*$, over the whole domain of ZF*-Ext., with the predicate “set” standing for $set$ as defined here which obviously is equivalent to $``\neg At"$), detailed as follows: \smallskip

\emph{Weak Extensionality:} $ \forall  X \forall  Y [set(X) \wedge set(Y)  \Rightarrow   ( \forall  Z(Z \in^*X \iff   Z \in^*Y)   \Rightarrow   X=^*Y)]$
Proof: Since X,Y are sets, so $Z \in^*X, Z \in^*Y$ becomes equivalent to $Z \in   X, Z \in   Y$ [Corollary 2], so we’ll have  $\forall  Z(Z \in   X \iff   Z \in   Y)$, which is $X=^*Y$ by definition.\smallskip

 $\in^*$-\emph{Foundation:}Proof is also straightforward from  $\in$   -Foundation. Since we have: \smallskip
 
 \noindent
 $\exists   x [ \exists   k(k \in^*x) \wedge \forall  y (y \in^*x   \Rightarrow   \exists   z (z \in^*y \wedge z  \in^* x))] \Rightarrow 
 \\\exists   x [ \exists   k(k \in   x) \wedge \forall  y (y \in  x   \Rightarrow   \exists   z (z \in y \wedge z \in x))]$ \smallskip
 
 (Corollary 2) ,then the violation of  $\in^*$-foundation will constitute a violation of  $\in$   -foundation.
\smallskip

 $\in^*$-\emph{Pairing:} for any objects $a,b$, by pairing we have a class $x$ such that $Pair(a,b,x)$, where the later formula is ($\forall  y (y \in    x  \iff    y=a \lor y=b)$), then by Proposition 1 we’ll have a set $x^*$ of all co-extensionals of $a$ and all co-extensionals of $b$.
\smallskip

 $\in^*$-\emph{Union:} for any set $x$, apply  $\in   $-Separation to get a set $x^*$ of all elements of $x$ that are sets, now take any $\in$-union of $x^*$ and we get an  $\in^*$-Union set of $x$. Now if $x$ is not a set, then the empty set will witness the existential quantifier for $\in^*$-union.
\smallskip

 $\in^*$-\emph{Power:} for any class $x$, take the Power class of $x$, and that would clearly be a $set$ [Subsidiary 2]
\smallskip

In all proofs of $\in^*$ versions of Pairing, Union and Power, we use Corollary 3 to finalize them \smallskip

$ \in^*$-\textbf{\emph{Replacement:}} Proof: Recall the axiom of $\in^*$-replacement would be:\smallskip

if $\phi(x,y)$ is a formula that only uses predicates $\in^*; =^*$ and not mentioning $B$, then each closure of: $$\forall A [\forall x \in^* A \exists!^* y (\phi(x,y)) \to \exists B: set(B) \land \forall y (y \in^* B \leftrightarrow \exists x \in^* A (\phi(x,y)))]$$; is an axiom. Where $``\exists !^*"$ is co-extensionality existential quantifier defined as: $$ \exists !^*y ( \phi) \iff \exists k \forall y (\phi(y) \iff y=^*k)$$ First we note that for any such formula $\phi$ that fulfills the stipulated antecedent, then by using Replacement* we can prove that: $$ \exists B: set(B) \land \forall y (y \in B \iff \exists x \in A \forall z (z \in y \iff \phi(x,z)))$$; the proof is pretty much straightforward: since we have $\exists   k \forall  y(\phi(x,z)\iff z=^*k)$, then we are replacing each $x \in A$ by a class $y$ of all objects $z$ that are co-extensional to some class $k$. By then each $y$ is a set! Now the existence of each set $y$ for each $x$ in $A$ is guaranteed by the condition that for each $x \in A$ each object $z$ that satisfies $\phi(x,z)$ must be co-extensional to some set $k$. Now from Pairing we have $\{k\}$ and from Proposition 1 we have the set of all co-extensionals of $k$ existing, and so there do exist at least one $y$ to replace each $x \in A$, so we get $B$ being nonempty if $A$ is nonempty, and the class $B$ is clearly a $set$ [Lemma 2]. Now from that we get the same formula above but with $y \in^* B$ (Corollary 2). What remains is to turn $x \in A$ to $x \in^* A$, now to see that this is the case, we'd observe that if $A$ is a set then the turning is in effect (Corollary 2), if $A$ is not a set, then $B$ would be $\in$-empty and so is a $set$, since the formula on the right of the biconditional is only using predicates $\in^*, =^*$, then we have $y \in^* B$ on the left (Corollary 3). QED
\smallskip

The proof of  $\in^*$-\textbf{Infinity} goes as follows: What is needed is to prove that: \smallskip

\noindent
$\exists N: \exists o \in^*N (\neg \exists z (z \in^* o) \land set(o))  \wedge \\\forall x \in^*N \exists y \in^* N  \forall z (z \in^* y \leftrightarrow z \in^*x \lor z=^* x)$; \smallskip

The idea is to define the usual Von Neumann ordinals (i.e. transitive classes of transitive classes such that no two distinct elements are co-extensional) as  $\in$- ordinals, and to define $ \in^*$-ordinals as transitive sets of transitive sets (Assuming Regularity), then to define a ``copying relation" between them such that for each  $\in$   - ordinal there is a copying relation $F$ that sends it to an  $\in^* $ - ordinal, and the latter shall be denoted as an  $\in^*$-ordinal copy of it. We’ll denote a usual finite $ \in$   -ordinal by $n$ and an  $\in^*$ - ordinal copy of it by $n^*$. Now we define the predicate ``ordinal copying relation" in the following manner: \smallskip

For all $F$, $F$ is an ordinal copying relation from $n$ to $m$ if and only if $F$ is a class of Kuratowski ordered pairs such that: \smallskip

\noindent
$\forall  a,b((a,b) \in   F\Rightarrow a \in   n \wedge b \in   m) \wedge   \\\forall  a,b,c,d( ( a,b) \in   F\wedge (c,d) \in   F \Rightarrow [a=c \iff b=^*d])$, \smallskip

that was the quasi-injectivity condition, also $F$ must satisfy the “quasi-surjectivity” condition which is: \smallskip

$[ \forall  a \in   n \exists   b \in   m( \exists   p \in   F(p=(a,b)))] \wedge [ \forall  b \in   m \exists   a \in   n( \exists   p \in   F(p=(a,b)))]$, \smallskip

therefore we'd say that $F$ is quasi-bijective. \smallskip

Now we also require ``quasi-isomorphism" condition, i.e. that of:
$$\forall  a,b,c,d((a,b) \in   F\wedge(c,d) \in   F\Rightarrow [a \in   c \iff b \in^*d])$$

 Now all  $\in^*$-empty ordinals are  $\in$   -empty sets, so a quasi-bijectivity ordinal copying relation trivially exists between them. So any first Von-Neumann ordinal will have an  $\in^*$-ordinal copy, which is simply itself! And actually any other empty set. Now we can prove by INDUCTION that for any “finite” Von-Neumann ordinal $n$, there is an $n^*$ copy of it, and also there is a set of all $n^*$ copies of it. This is easy since this is the case for every empty set $n$, also for any $n$ that has an $n^*$ and a set $F$ that is an ordinal copying relation between it and that $n^*$, then clearly $n+1$ will have an $(n+1)^*$ and a set $G$ that is an ordinal copying relation between them, since $G$ will be a union of $F$ with a Cartesian product of an $ \{n\}$ and an  $\{x|x=^*n^*\}$, and $(n+1)^*$ would be the range of $G$. That induction holds here is due to “Regularity”, since the von Neumann ordinals are all  $\in$   -well founded, so no descending membership chain is allowed, so induction holds in the usual manner.
\smallskip

So in nutshell for every  $\in$-ordinal $n$ that has an $\in^{*}$-ordinal copy $n^{*}$, each successor $n+1$ of $n$ will be easily shown to be sent by some ordinal copying relation to each successor $(n+1)^{*}$ of $n^{*}$, where $(n+1)^{*}$ is a union set of $n^{*}$ and a set $\{x|x=^{*} n^{*}\}$. Now for any set $n^{*}$ we have any set  $\{x|x=^*n^{*}\}$ being a set that contains ALL co-extensionals of any set $n^{*}$ because simply ALL $n^{*}$ sets are co-extensional (Induction)! So this proves that for every finite von Neumann ordinal $n$, there is a set of all  $\in^{*}$-ordinal copies of it.
\smallskip

Now by Replacement* we can send each element $n$ of the set of all finite Von Neumann ordinals to all sets $y$ of all of $n^*$ copies. Now a union of this replacement set would be a set that corresponds to an  $\in^*$-$\omega$. Or another way is to use Replacement* to replace each element $n$ of the set of all finite Von Neumann ordinals by each set $y$ of all  $\in^*$-ordinal copies of elements of $n$, in other words each $n$ would by replaced by all $n^*$ sets, and so the resulting set from Replacement* would be directly an $\in^*$-$\omega$ set. QED
\smallskip

So far we proved that ZF*-Extensionality would interpret ZFA. And this is known to interpret full ZF. However, we’ll give a line of interpretation along that method:
Now to prove full ZF we need to restrict ourselves to the domain of hereditarily $sets$, or more appropriately termed as “Pure sets”, i.e. classes that are $sets$ and having transitive closures in which every element is a $set$. It would be seen that this theory would prove all sentences of ZF with its variables so restricted and its predicates being  $\in^*, =^*$ instead of epsilon and equality.
\smallskip

We need to prove the existence of a transitive closure for each class. To do that we define the notion $\bigcup(x,n,y)$ to mean $y$ is an $n^{th}$ union of $x$, where $n$ is a finite Von Neumann ordinal. Now this is defined as:
\smallskip

$\bigcup(x,n,y)  \iff \exists k(k \in n)  \wedge \\
 \exists   set\ F \ ( \forall  p \in   F  \exists   a  \in    n  \exists   b (p=(a,b)) 
\wedge \\  \forall  m \in   n  \exists   p  \in    F  \exists   q (p=(m,q))
\wedge \\  \forall  a,b,c,d ((a,b)  \in    F \wedge  (c,d)  \in    F   \Rightarrow   [S(a,c)  \iff    \bigcup(b,d)])
\wedge \\ \exists   p  \in    F  (p=(first(n),x)) 
\wedge \\  \exists   p  \in    F (p=(last(n),y))$ \smallskip

Where $S(a,c)$ stands for $c$ is a successor of $a$ (i.e.; $c$ is an $ a \cup \{a\}$ class), $\bigcup(b,d)$ signify $d$ is a class union of $b$.\smallskip

Replace (through Replacement*) from the set of all non-empty finite von Neumann’s $ N^+$ where each $n$  $\in$    $N^+$ is sent to each $y$ such that $\bigcup(x,n,y)$. This can be done since for a certain $x$, all $y$'s for each $n$ are co-extensional(Induction)! Then we take the union of the resulting set and we get a class $t$ that is a transitive closure class of $x$, this is denoted as $TC(x,t)$. Now we define: $$ x \ is \ hereditarily \ set   \iff    x \ is \ a \ set \wedge  \exists  t (TC(x,t)   \wedge   \forall  y  \in    t(y \ is \ a \ set))$$ \textbf{Proposition 2:}  $\forall  x (set(x) \wedge  \forall  y (y  \in    x   \Rightarrow   y \ is \ a \  pure \ set) \Rightarrow x \ is \ a \  pure \ set)$ \smallskip                            
Proof: every element of any transitive closure of $x$ is either an element of $x$ or an element of a transitive closure of an element of $x$; since all of those are sets and $x$ itself is a set, then $x$ is a Pure set.
\smallskip

Now of the above Axioms interpreting ZFA, all of them when their parameters are pure sets, then the asserted to exist sets are pure! Because each output $y$ of all their defining formulas is provably a pure set, so the asserted set of all these outputs would be accordingly pure (Proposition2). Also the pure set realm is transitive (i.e.; an element of a pure set is a pure set)! So the whole pure set world is shown to be closed under these Axioms. And of course full Extensionality is directly satisfied! Thus interpreting full ZF.
\smallskip

Now what needs to be get rid of, is the axiom of Regularity, this is done by defining a Hierarchy building relation in manner similar to relation $F$ above but using Powers instead of Unions, and stipulating that each ``limit” ordinal in a domain of $F$ is sent by $F$ to a union of all $F$-images  of members of the domain strictly lower than that limit, and also we must stipulate that $F$ sends the first ordinal in its domain to an empty set. Now a well founded set can then be defined as an element of an element in the range of a hierarchy building function. Then induction can be used over the indices of those power stages since they are ordinals which are \emph{built-in} well founded classes (i.e. by adding the requirement of being $\in$-well founded to their above mentioned definition), and accordingly the property of regularity would be proved over all sets in the stages of these hierarchies, and so induction can be used, and the whole proof goes through within those cumulative hierarchies in the usual manner \smallskip

\textbf{An aside:} I think the simplest modification on the conventional axiom schema of Replacement that can interpret the above scheme is the following: $$[\forall X \exists Z \forall Y (\phi(X,Y) \Rightarrow Y \in Z)] \Rightarrow \forall A \exists B \forall Y (Y \in B \iff \exists X \in A (\phi(X,Y)))$$ This can be labeled as ``Union schema". This would prove (with Power) Pairing, Union, Separation, and Replacement*. Now if we change $Y\in Z$ in the pre-conditional to $Y \subset Z$ then together with Infinity we get a theory that would interpret ZF over the well founded set realm of it.

\section{ A Separate proof}

Azrield Levy had reported an important result due to Dana Scott. Here is a quote from his artcile:\smallskip

``If ZF* means ZF but with replacement written as follows, then ZF*-Ext. interprets ZF. \smallskip

Replacement*: $$\forall A[\forall x,y,z (\phi(x,y) \wedge \phi(x,z) \Rightarrow y=^*z) \Rightarrow \exists B \forall y (y \in B \iff \exists x \in A \phi(x,y))]$$ where $ y=^*z$ is defined as [$\forall m (m \in y \iff m \in z)$]"\smallskip

This is a theorem here! Just let $\phi(x,z)$ in Replacement* be $ \exists k (\phi(x,k) \wedge z \in k) $, then $\in$-separate on the resulting set by  $\exists x \in A (\phi(x,y))$. 

To be noted is that all proofs presented here can depend on this principle without the need consult our version  Replacement*.

\section{References}

\noindent
[1] Al-Johar, Z; Holmes M.R., Acyclic Comprehension is equal to Stratified Comprehension, Preprint 2011. http://zaljohar.tripod.com/acycliccomp.pdf

\noindent
[2] Al-Johar, Zuhair; Holmes, M. Randall; and Bowler, Nathan. (2014). "The Axiom Scheme of Acyclic Comprehension". Notre Dame Journal of Formal Logic, 55(1), 11-24. http://dx.doi.org/10.1215/00294527-2377851

\noindent
[3] Crabb`e, Marcel, “On NFU”. Notre Dame Journal of Formal Logic 33 (1992), pp 112-119.

\noindent
[4] R. O. Gandy, On the axiom of extensionality I. J. Symbolic Logic 21, 36 – 48 (1956).

\noindent
[5] R. O. Gandy, On the axiom of extensionality II. J. Symbolic Logic 24, 287 – 300 (1959).

\noindent
[6] R. Hinnion, Extensional quotients of structures and applications to the study of the axiom of extensionality. Bull. Soc. Math. Belgique, S`erie B, 33 (1981).

\noindent
[7] D. Scott, More on the axiom of extensionality. In: Essays on the Foundations of Mathematics, 
pp. 115 – 131 (Magness Press, Jerusalem 1961). 

\noindent
[8] Azriel Levy, R. O. Gandy. On the axiom of extensionality. The journal of symbolic logic, vol. 21 (1956), pp. 36–48, and vol. 24 no. 4 (for 1959, pub. 1961), pp. 287–300. - Dana Scott. More on the axiom of extensionality. Essays on the foundations of mathematics, dedicated to A. A. Fraenkel on his seventieth anniversary, edited by Y. Bar-Hillel, E. I. J. Poznanski, M. O. Rabin, and A. Robinson for The Hebrew University of Jerusalem, Magnes Press, Jerusalem 1961, and North-Holland Publishing Company, Amsterdam1962, pp. 115–131. J. Symbolic Logic 29(3), p.142, (1964)

\end{document}